\documentclass[12pt]{article}%
\usepackage{amsmath}
\usepackage{url}
\usepackage{amsfonts}
\usepackage{amssymb}
\usepackage{graphicx}
\usepackage{color}
\usepackage{url}
\usepackage{hyperref}%
\providecommand{\U}[1]{\protect\rule{.1in}{.1in}}
\providecommand{\U}[1]{\protect\rule{.1in}{.1in}}

\begin{document}

\title{Evidence-based teaching: how do we all get there?}
\author{David Pengelley\\Oregon State University, Corvallis, OR 97331\\davidp@nmsu.edu
\and Dev Sinha\\University of Oregon, Eugene, OR 97403\\dps@uoregon.edu}
\date{June 3, 2019}
\maketitle

There are compelling reasons to shift our pedagogy toward evidence-based
active learning methods that substantially improve student success, and now
plenty of resources to aid in that shift. These include the recent
\href{https://www.cbmsweb.org/cbms-position-statements/}{CBMS Statement on
Active Learning},
\href{https://www.maa.org/programs-and-communities/curriculum%20resources/instructional-practices-guide}{MAA
Instructional Practices Guide} (IPG), and
\href{http://math.mit.edu/seminars/emes/index.html}{MIT Electronic Seminar on
Mathematics Education} (see also Resources below). But implementation is
neither quick nor easy. There are still plenty of individual, institutional,
cultural, and professional obstacles, along with wonderful opportunities.

At the 2019 Joint Mathematics Meetings we co-organized a guided discussion --
an ``un-panel'' -- sponsored by the American Mathematical Society's\ Committee
on Education in order to stimulate the process of our community moving toward
active learning in our teaching pedagogy. Seventy participants with fifteen
discussion leaders expanded an initial list of issues, and considered
questions around both challenges and opportunities. Here we summarize from
these discussions, suggesting areas for collaborative efforts ranging from
local colleagues and educational institutions to national and global
professional societies.

Fourteen issues were identified, and although there is overlap between them,
we summarize participant comments separately for each, mostly direct quotes
from notes taken by discussion leaders (more quotes are available in our
report\footnote{\href{https://web.nmsu.edu/~davidp/how-do-we-all-get-there-followup.pdf}{https://web.nmsu.edu/\symbol{126}%
davidp/how-do-we-all-get-there-followup.pdf}}):

\begin{enumerate}
\item Training graduate students and early career mathematicians: Feelings in
a word: \textquotedblleft Isolated.\textquotedblright\ There is a need for
mentors, encouragement, community, local discussion, concrete practices.
\textquotedblleft Need for a ramp on - how do you start learning to do this
and how do you feel supported in the learning?\textquotedblright

\item Developing departmental experts who can lead and mentor: ``More people
need to rise as leaders - not be appointed necessarily - but gather groups
together to talk about and try active learning.'' There is a ``lack of
expertise, a lack of confidence, and a lack of time. If there were a visible
group of people doing active learning at the institution, a person could be
mentored,'' but without such a group we won't see progress.

\item Offering wide-scale programming for department chairs: Making
\textquotedblleft recommendations to administrators giving active learning
pedagogies the imprimatur of the mathematics professional organizations (in
addition to the marvelous CBMS statement on active learning) would help
nationally.\textquotedblright\ \textquotedblleft Engage the department faculty
in discussions on the case for adopting active learning pedagogy and the
consequent need for the faculty who are not currently using it to give it
consideration.\textquotedblright\ Have \textquotedblleft chairs/leadership
encourage an MAA IPG-reading club.\textquotedblright\ \textquotedblleft
Highlight active learning pedagogy as the gold standard for effective methods
of instruction.\textquotedblright

\item Updating inventory tools of teaching practices for observations and
training: \textquotedblleft We need an inventory tool building on the MAA IPG.
This is an area where our professional societies can make a big contribution
by supplying this for the community.\textquotedblright

\item Shifting program evaluation towards active learning and deeper, more
authentic learning outcomes: \textquotedblleft Active learning goes
hand-in-hand with a shift from low-level recall expectations of students to
more valuable outcomes. We should be going from assessments which just test
whether one can recall and perform appropriate techniques of integration to
instead asking for some proofs (even at entry level) or authentic applications
(given only a context one has to see how an integral would apply, discuss
assumptions and accuracy, etc.) or even say historical development. If one
shifts attention to these kinds of outcomes, then classrooms obviously need to
put more emphasis on students practicing these, which means a shift to active
learning.\textquotedblright

\item Large lectures, and the challenges they present for interaction,
including individual feedback, group work, and whole-class discussion: A good,
proven strategy is having \textquotedblleft undergraduate learning assistants
(see Resources below) and graduate teaching assistants in the active learning
classroom.\textquotedblright

\item Dissemination of teaching materials currently biased towards lecture
format: We should \textquotedblleft gather/curate materials that have been
developed for specific courses, but are not in wide circulation, to help
minimize the amount of work and effort required to adopt and implement active
learning pedagogies.\textquotedblright\ Or think more broadly about
\textquotedblleft new publishing models and systems for pedagogical materials
including worksheets, think-pair-share questions, applets, discussion of areas
of difficulty and common student responses and plans for addressing them, and
other materials which support active learning, along with textbooks, which
combine features of the arXiv, MathReviews, MathOverflow, Curated Courses \&
UTMOST, Webwork, CoMInDS training materials site. Financed by course fees,
with AMS/MAA as fiscal intermediary (new financial model), keeping more money
in the community and saving students money. In the meantime... perhaps have
more robust, high-profile blog networks (\#MTBoS) where people share key
experiences and favorite materials and approaches.\textquotedblright

\item Culture, inertia, and incentives: Institutions and the community should
\textquotedblleft create incentives for trying new pedagogical strategies
(including \$\$ for professional development).\textquotedblright\ We should
have \textquotedblleft discussions about changing the reward system so that
the faculty who adopt active learning pedagogies get credit for doing
so.\textquotedblright\ There should be a \textquotedblleft shift in the
institutional culture from one in which the traditional lecture method is
accepted as the default/standard method of instruction to one in which active
learning methods become the preferred method of instruction.\textquotedblright%
\ \textquotedblleft Internal grants for research can be used to model novel
programs that include internal grants for teaching skills
development.\textquotedblright\ \textquotedblleft We have lots of challenges,
including many structural ones (physical spaces, cultural inertia, incentive
structures, lack of professional development built into the higher education
model, large number of adjunct and term instructor positions) - much of this
is, at its core, really about the R1 business model.\textquotedblright%
\ \textquotedblleft Institutional change is essential for wide scale adoption
of active learning pedagogies.\textquotedblright\ We should be
\textquotedblleft creating department-level statements (like the CBMS
statement) about active learning.\textquotedblright\ \textquotedblleft Change
the standards used for evaluating faculty to place greater value on active
learning.\textquotedblright\ Institutions should adopt \textquotedblleft
student evaluations that address active learning, and are less biased against
it.\textquotedblright\ \textquotedblleft For teaching evaluation:
\href{https://salgsite.net/}{SALG} and
\href{https://www.ideaedu.org/Services/Services-to-Improve-Teaching-and-Learning/Student-Ratings-of-Instruction}{IDEA}%
.\textquotedblright\ \textquotedblleft At conferences, MAA/AMS should have
pedagogy sessions/component tied to each math course area. Don't schedule
against the research sessions.\textquotedblright\ \textquotedblleft When we go
to give a talk and work on research, ask to sit in on
classes.\textquotedblright\ We should be \textquotedblleft normalizing
education talks at seminars/colloquia; lists of visiting lecturers to give
these talks.\textquotedblright\ We should create \textquotedblleft small,
working communities of practice in both research and teaching. On the research
side, looks like Women In... model of collaborative work on an expert-chosen
problem and/or REUF. But also have pedagogy leader(s) to focus on producing
some materials (worksheets, texts, applets, discussions of practice)
addressing a problem of pedagogical practice as well. In the meantime ...
let's make it a common practice to have a pedagogical session as part of
research conferences, especially graduate training
conferences.\textquotedblright

\item Facilities: Need facilities that are compatible with \textquotedblleft
group work and board work (group tables, boards on walls).\textquotedblright

\item Informed support as resource: \textquotedblleft Offer/provide
professional development to any faculty member who is interested in adopting
active learning instructional methods in their classes.\textquotedblright%
\ \textquotedblleft The institutional center for excellence in teaching and
learning is a good resource for active learning pedagogies.\textquotedblright%
\ \textquotedblleft Academy of Inquiry Based Learning (IBL) website is a
terrific resource. The IBL Video Series are phenomenal and showcase models for
novices to use. The site also features textbooks, links, and course materials
developed by IBL experts.\textquotedblright\ \textquotedblleft A teaching
seminar for faculty?\textquotedblright\ \textquotedblleft Focus more on what
$\emph{students}$ do, rather than faculty.\textquotedblright%
\ \textquotedblleft Information on which universities have changed their
reward structure to allow active learning pedagogies to count as a valued
activity for purposes of tenure and promotion.\textquotedblright%
\ \textquotedblleft A `kit' that goes to departments to show how to get things
started.\textquotedblright\ \textquotedblleft A department seminar on
educational practices; AMS Committee on Education could create a list of
expert speakers who would be willing to come speak at such seminars, lending
professional credence to issues of teaching.\textquotedblright%
\ \textquotedblleft RUME\ folks could do more interpreting of math ed research
for mathematicians, e.g., what does this mean for me when I go into my
calculus class?\textquotedblright\ \textquotedblleft Formalizing via workshops
and conferences, book club groups to learn more (e.g., MAA IPG), MAA section
meetings, NExT sections, virtual resources (MIT seminar).\textquotedblright

\item Finding resources and the time involved: \textquotedblleft Start small;
identify colleagues; find allies.\textquotedblright\ \textquotedblleft Use a
Reading Group to focus on the MAA IP Guide.\textquotedblright\ Cultivate
\textquotedblleft institutional support for class setup.\textquotedblright%
\ \textquotedblleft People need to see lots of examples of what can be
done.\textquotedblright

\item Collaborating with peers, seeing classrooms: \textquotedblleft Give
common tests (write together), collaborate on grants, attend IBL workshops
together, start with those who are open and willing, ask to come sit in their
class, start with just 10\%, start with a colleague's materials, then tweak
and make it your own.\textquotedblright\ \textquotedblleft Critical mass of
folks interested in interactive, engaged teaching.\textquotedblright%
\ Unfortunately, \textquotedblleft we are not open about what we do in our own
classes.\textquotedblright\ \textquotedblleft Why do people refuse to be
observed by you? Fear of being judged?\textquotedblright\ \textquotedblleft Be
more open about sharing what goes on in the classroom.\textquotedblright%
\ \textquotedblleft Sit in on each other's classes.\textquotedblright%
\ \textquotedblleft Collaborate on making teaching public.\textquotedblright%
\ \textquotedblleft Create incentives to making teaching
public.\textquotedblright

\item What to do to bring in colleagues?: \textquotedblleft It isn't all or
nothing.\textquotedblright\ \textquotedblleft Talk about what students do and
need; this is less threatening.\textquotedblright\ \textquotedblleft Compile
and share data showing that active learning pedagogy is more effective than
other pedagogies. The more local the data the better.\textquotedblright%
\ \textquotedblleft Pure mathematicians are sometimes not convinced by
mathematics education research.\textquotedblright

\item Leading from below, e.g., grad student or faculty member:
\textquotedblleft Support leading from below (grad students, postdocs, early
career).\textquotedblright\ Think about \textquotedblleft CV
building.\textquotedblright\ \textquotedblleft Lack of influence (grad
students/postdocs).\textquotedblright\ The general feeling is that one is
\textquotedblleft not expected/invited to contribute as a grad.
student.\textquotedblright\ Still predominantly hear \textquotedblleft advice
from mentors that too much emphasis on teaching makes one look like a less
valuable researcher.\textquotedblright\ \textquotedblleft Being valued, being
part of the conversation.\textquotedblright
\end{enumerate}

Clearly there are many factors influencing pedagogical change, and many scales
at which they may occur. Change can happen at the wholly individual level, or
as part of an institutional team, or through collaboration with like-minded
colleagues at other institutions, or at the broader scale in which our
professional societies should exert key leadership. We end with a point made
recently to us by Uri Treisman, that changing structures is a key step in
widespread reform. The participant discussions quoted above reveal that there
are many structures for us to change at the departmental, institutional, and
more global levels in order to achieve the goal of evidence-based active
learning pedagogy in our teaching.

We encourage everyone to continue such conversations and find their next step
in this process, and to share their experience with the authors of this
article so that we may further disseminate progress.

\bigskip

\begin{center}
{\large \textbf{Resources} }
\end{center}

\begin{itemize}
\item NSF press release
\emph{\href{https://www.nsf.gov/news/news_summ.jsp?cntn_id=131403}{\emph{Enough
with the lecturing}}}\ (2014)

\item CBMS statement
\emph{\href{https://www.cbmsweb.org/cbms-position-statements/}{\emph{Active
Learning in Post-Secondary Mathematics Education}}\ }(2016)

\item National Research Council:
\emph{\href{https://www.nap.edu/catalog/18687/reaching-students-what-research-says-about-effective-instruction-in-undergraduate}{\emph{Reaching
Students: What Research Says About Effective Instruction in Undergraduate
Science and Engineering}}} (2015)

\item
\href{https://www.maa.org/programs-and-communities/curriculum%20resources/instructional-practices-guide}{MAA
Instructional Practices Guide} (2017)

\item
\href{https://www.amstat.org/asa/education/Guidelines-for-Assessment-and-Instruction-in-Statistics-Education-Reports.aspx}{Guidelines
for Assessment and Instruction in Statistics Education} (GAISE) (2016)

\item \href{http://math.mit.edu/seminars/emes/index.html}{MIT Electronic
Seminar on Mathematics Education}

\item \href{https://blogs.ams.org/matheducation/}{AMS Blog on Teaching and
Learning Mathematics}

\item \href{http://cominds.maa.org/}{College Mathematics Instructor
Development Source (CoMInDS)}, MAA

\item
\href{https://www.maa.org/programs-and-communities/professional-development/project-next}{Project
NExT (New Experiences in Teaching)}, MAA; \& Project ACCCESS, AMATYC

\item
\href{https://www.aau.edu/education-service/undergraduate-education/undergraduate-stem-education-initiative}{Undergraduate
STEM Education Initiative}, AAU

\item
\href{http://www.aplu.org/projects-and-initiatives/stem-education/seminal/index.html}{Student
Engagement in Mathematics through an Institutional Network for Active
Learning} (SEMINAL), APLU

\item \href{http://www.inquirybasedlearning.org/}{Academy of Inquiry Based
Learning}

\item \href{http://www.jiblm.org/}{Journal of Inquiry-Based Learning in
Mathematics}

\item \href{http://www.mathlearningbyinquiry.org/}{Initiative for Mathematics
Learning by Inquiry} (MLI)

\item \href{http://sigmaa.maa.org/ibl/}{IBL Special Interest Group of the MAA}

\item \href{https://www.tpsemath.org/}{Transforming Post-Secondary Education
in Mathematics} (TPSE Math)

\item \href{https://www.pogil.org/}{Process Oriented Guided Inquiry Learning} (POGIL)

\item
\href{https://www.maa.org/programs-and-communities/member-communities-innovative-teaching-exchange}{Innovative
Teaching Exchange}, MAA

\item
\href{http://www.ams.org/programs/edu-support/Innovations-College-Level}{Innovations
in College-Level Mathematics Teaching}, AMS

\item \href{https://www.learningassistantalliance.org/}{Learning Assistant
Alliance}; \&
\href{https://www.colorado.edu/program/learningassistant/}{Learning Assistant
Program} (UC Boulder)

\item \href{https://ascnhighered.org/index.html}{Accelerating Systemic Change
Network} (ASCN)

\item \href{https://www.nsf.gov/div/index.jsp?div=DUE}{NSF Division of
Undergraduate Education}

\item
\href{https://www.maa.org/programs-and-communities/curriculum%20resources/progress-through-calculus}{MAA
Progress Through Calculus studies}

\item \href{https://amatyc.site-ym.com/page/Webinars}{AMATYC Webinar series}

\item \href{https://www.cirtl.net/}{Center for the Integration of Research,
Teaching, and Learning}

\item Active learning resources from David Pengelley:
\texttt{\href{https://web.nmsu.edu/~davidp/}{\texttt{https://web.nmsu.edu/\symbol{126}%
davidp/}}}

\item Active learning resources from Robin Pemantle, including pedagogical
tips, and materials for calculus and teacher preparation courses:

\texttt{\href{https://www.math.upenn.edu/~pemantle/Active-resources.html}{https://www.math.upenn.edu/\symbol{126}%
pemantle/Active-resources.html}}

\item Active learning course materials from Matt Boelkins, Steve Schlicker,
and Ted Sundstrom, for courses in

Active Calculus (single and multivariable):
\texttt{\href{https://activecalculus.org/}{https://activecalculus.org/}}

Active Preparation for Calculus:
\texttt{\href{https://gvsu.edu/s/0Ui}{https://gvsu.edu/s/0Ui}}

Mathematical Reasoning: Writing and Proof:
\texttt{\href{https://scholarworks.gvsu.edu/books/9/}{https://scholarworks.gvsu.edu/books/9/}%
}
\end{itemize}

\end{document}